\documentclass[12pt]{article}
\usepackage{amsmath}
\usepackage{latexsym}
\usepackage{amssymb}
\usepackage{a4wide}
\newtheorem{thm}{Theorem}[section]
\newtheorem{la}[thm]{Lemma}
\newtheorem{Defn}[thm]{Definition}
\newtheorem{Remark}[thm]{Remark}
\newtheorem{Note}[thm]{Note}
\newtheorem{prop}[thm]{Proposition}

\newtheorem{cor}[thm]{Corollary}
\newtheorem{Example}[thm]{Example}
\newtheorem{Examples}[thm]{Examples}
\newtheorem{Problems}[thm]{Problems}

\newtheorem{Problem}[thm]{Problem}
\newtheorem{Number}[thm]{\!\!}
\newenvironment{defn}{\begin{Defn}\rm}{\end{Defn}}

\newenvironment{example}{\begin{Example}\rm}{\end{Example}}

\newenvironment{rem}{\begin{Remark}\rm}{\end{Remark}}
\newenvironment{numba}{\begin{Number}\rm}{\end{Number}}
\newenvironment{proof}{{\noindent\bf Proof.}}%
                  {\nopagebreak\hspace*{\fill}$\Box$\medskip\medskip\par}   
\newcommand{\Punkt}{\nopagebreak\hspace*{\fill}$\Box$}

\newcommand{\wb}{\overline}
\newcommand{\ve}{\varepsilon}

\newcommand{\tensor}{\otimes}
\newcommand{\impl}{\Rightarrow}

\newcommand{\mto}{\mapsto}

\newcommand{\isom}{\cong}

\newcommand{\N}{{\mathbb N}}
\newcommand{\R}{{\mathbb R}}

\newcommand{\K}{{\mathbb K}}

\newcommand{\bD}{{\mathbb D}}

\newcommand{\C}{{\mathbb C}}

\DeclareMathOperator{\ad}{ad}

\newcommand{\cO}{{\mathcal O}}

\newcommand{\cS}{{\mathcal S}}
\newcommand{\cg}{{\mathfrak g}}

\newcommand{\dl}{{\displaystyle \lim_{\longrightarrow}}}
\newcommand{\pl}{{\displaystyle \lim_{\longleftarrow}}}

\newcommand{\sub}{\subseteq}

\DeclareMathOperator{\id}{id}

\newcommand{\cT}{{\mathcal T}}

\DeclareMathOperator{\Spann}{span}

\DeclareMathOperator{\eval}{eval}
\newcommand{\sbull}{{\scriptscriptstyle \bullet}}

\newcommand{\aeq}{\Leftrightarrow}

\begin{document}
\renewcommand{\thefootnote}{\fnsymbol{footnote}}
\begin{center}
{\Large\bf
Applications of hypocontinuous bilinear maps\\[1.6mm]
in infinite-dimensional
differential calculus}\\[4.6mm]
{\bf Helge Gl\"{o}ckner\footnote{This article is partially
based on the unpublished preprint~\cite{FAM},
the preparation of which was supported
by the German Research Foundation (DFG, FOR 363/1-1).}}\vspace{3mm}
\end{center}
\renewcommand{\thefootnote}{\arabic{footnote}}
\setcounter{footnote}{0}
\begin{abstract}\vspace{1.4mm}
\hspace*{-7.2 mm}
Paradigms of
bilinear maps $\beta\colon E_1\times E_2\to F$
between locally convex spaces
(like evaluation or composition)
are not continuous,
but merely hypocontinuous.
We describe situations
where, nonetheless,
compositions of $\beta$ with
Keller $C^n_c$-maps (on suitable domains)
are~$C^n_c$.
Our main applications
concern holomorphic families of operators,
and the foundations of locally convex Poisson
vector spaces.\vspace{2mm}
\end{abstract}
{\footnotesize {\em Classification}:
26E15, 
26E20 (primary); 
17B63, 
22E65, 
46A32, 
46G20, 
46T25 (secondary).\\[1.5mm] 
}\vspace{3mm}
\begin{center}
{\Large\bf Introduction}\vspace{.6mm}
\end{center}
If $\beta\colon E_1\times E_2\to F$
is a continuous bilinear map between
locally convex spaces, then~$\beta$ is smooth
and hence $\beta\circ f\colon U\to F$
is smooth for each smooth map
$f\colon U\to E_1\times E_2$
on an open subset~$U$ of a locally convex space.\\[3mm]
Unfortunately, many
bilinear mappings of interest
are discontinuous.
For example, it is known
that the evaluation map
$E'\times E\to \R$,
$(\lambda,x)\mto \lambda(x)$
is discontinuous for each
locally convex vector topology
on~$E'$, if~$E$ is a non-normable
locally convex space (cf.\ \cite{Mai}).
Hence also the composition map
$L(F,G) \times L(E,F) \to L(E,G)$,
$(A,B )\mto A \circ B$
is discontinuous,
for any non-normable
locally convex space~$F$,
locally convex spaces
$E,G\not=\{0\}$,
and any locally
convex vector topologies on
$L(F,G)$, $L(E,F)$ and
$L(E,G)$ such that the maps
$F\to L(E,F)$, $y\mto y\tensor \lambda$
and $L(E,G)\to G$, $A\mto A(x)$
are continuous for some
$\lambda\in E'$
and some $x\in E$ with $\lambda(x)\not=0$,
where $(y\tensor\lambda)(z):= \lambda(z)y$
(see Remark~\ref{nachtrg}).\\[3mm]
Nonetheless, both evaluation and composition
do show a certain weakened
continuity property, namely \emph{hypocontinuity}.
So far, hypocontinuity arguments
have been
used in differential calculus
on Fr\'{e}chet spaces in some isolated cases
(cf.\ \cite{FAM}, \cite{HYP} and \cite{Wur}).
In this article, we distill
a simple, but useful general principle
from these arguments
(which is a variant of a result from~\cite{Tho}).
Let us call a Hausdorff
topological space~$X$ a \emph{$k^\infty$-space}
if $X^n$ is a $k$-space for each $n\in \N$.
Our observation
(recorded in Theorem~\ref{mainobs})
is the following:\\[3mm]
\emph{If a bilinear map $\beta\colon E_1\times E_2\to F$
is hypocontinuous with respect to
compact subsets of $E_1$ or~$E_2$
and $f\colon U\to E_1\times E_2$
is a $C^n$-map on an open subset~$U$ of a locally convex space~$X$
which is a $k^\infty$-space,
then $\beta\circ f\colon U\to F$ is
a $C^n$-map.}\vspace{1.3mm}\pagebreak

\noindent
As a byproduct, we obtain an affirmative solution
to an old open problem by Serge
Lang (see Corollary~\ref{corlang}).
Our main applications
concern two areas.\\[3mm]
{\bf Application 1: Holomorphic families of operators.}\\[1.3mm]
In Section~\ref{secholfam},
we apply our results
to holomorphic
families of operators,
i.e., holomorphic maps
$U\to L(E,F)$ on an open set~$U\sub \C$
(or $U\sub X$ for a suitable
locally convex space~$X$).
We obtain generalizations
(and simpler proofs)
for various results
formulated in~\cite{BaO}.\\[3mm]
{\bf Application 2:
Locally convex Poisson vector spaces.}\\[1.3mm]
Finite-dimensional
Poisson vector spaces
are encountered naturally in finite-dimensional
Lie theory as the dual spaces $\cg^*$
of finite-dimensional Lie algebras.
They give rise to a distribution on~$\cg^*$
whose maximal integral manifolds
are the coadjoint orbits
of the corresponding simply connected
Lie group~$G$. These are known
to play an important role
in the representation
theory of~$G$ (by Kirillov's
orbit philosophy).\\[2.5mm]
The study of infinite-dimensional
Poisson vector spaces
(and Poisson manifolds)
only began recently with
works of Odzijewicz and
Ratiu concerning the Banach case
(see \cite{OR1}, \cite{OR2}).
On Neeb's initiative,
a setting of
locally convex Poisson vector spaces
(which need not be Banach spaces)
has recently been developed (\cite{GLN}, \cite{Ne2}).
In Section~\ref{secpoisson},
we preview this framework
and explain how
hypocontinuity arguments
can be used to overcome the
analytic problems arising beyond the Banach case.
In particular, hypocontinuity
is the crucial tool
needed to define
the Poisson bracket and Hamiltonian vector
fields.
In Section~\ref{nwsec1},
we describe situations where
the Poisson bracket is hypocontinuous
or even continuous.
In the most relevant case,
we also prove continuity
of the linear map which takes
a smooth function
to the corresponding Hamiltonian
vector field (Section~\ref{nwsec2}).\\[2.5mm]
Having thus secured the foundations,
natural next steps will be the investigation
of co\-adjoint orbits
for examples of non-Banach,
infinite-dimensional Lie groups
and geometric quantization in this context.
These will be undertaken in \cite{GLN},
\cite{Ne2} and later research.\\[2.8mm]
Let us remark that
an alternative path leading
to Theorem~\ref{mainobs}
is to prove,
in a first step, smoothness of hypocontinuous
bilinear maps (with respect to compact sets)
in an alternative
sense, replacing ``continuity''
by ``continuity on each compact subset''
in the definition of a $C^n$-map
(see~\cite[Theorem~4.1]{Tho}).\footnote{Compare
also \cite{Sei} for the use
of Kelleyfications in differential calculus.}
The second step
is to observe that such $C^n$-maps
coincide with ordinary $C^n$-maps
(i.e., Keller $C^n_c$-maps)
if the domain~is~a~\mbox{$k^\infty$-space.}\\[2.8mm]
\emph{Acknowledgement.}
The referee
of~\cite{FAM} suggested some valuable
references to the literature.
%
%
%
%
%
%
%
%
\section{Preliminaries and basic facts}\label{scprelim}
Throughout this article,
$\K\in \{\R,\C\}$
and $\bD:=\{z\in \K\colon
|z|\leq 1\}$.
As the default,
the letters
$E$, $E_1$, $E_2$,
$F$ and~$G$ denote
locally convex topological
$\K$-vector spaces.
When speaking of linear
or bilinear maps,
we mean $\K$-linear
(resp., $\K$-bilinear) maps.
A subset $U\sub E$ is called
\emph{balanced} if $\bD U\sub U$.
The locally convex spaces
considered need not be Hausdorff,
but whenever they serve as the domain
or range of a differentiable
map, we tacitly assume the Hausdorff
property.
We are working in a setting
of infinite-dimensional differential
calculus known
as Keller's $C^n_c$-theory
(see, e.g., \cite{RES} or \cite{GaN}
for streamlined introductions).
\begin{defn}\label{C11}
Let $E$ and~$F$ be locally convex
spaces over
$\K\in \{\R,\C\}$,
$U\sub E$ be open
and $f\colon U\to F$ be a map.
We say that $f$ is $C^0_\K$ if $f$ is continuous.
The map~$f$ is called $C^1_\K$ if it is continuous,
the limit
\[
df(x,y)\;=\;\lim_{t\to 0}\frac{f(x+ty)-f(x)}{t}
\]
exists for all $x\in U$
and all $y\in E$
(with $0\not=t\in \K$ sufficiently small),
and the map $df\colon U\times E\to F$
is continuous.
Given $n\in \N$, we say that
$f$ is $C^{n+1}_\K$ if $f$ is
$C^1_\K$ and $df\colon U\times E\to F$
is~$C^n_\K$.
We say that~$f$ is $C^\infty_\K$
if~$f$ is~$C^n_\K$ for each $n\in \N_0$.
If~$\K$ is understood,
we simply write $C^n$ instead of $C^n_\K$,
for $n\in \N_0\cup\{\infty\}$.
\end{defn}
If $f\colon E\supseteq U\to F$
is $C^1_\K$, then $f'(x):=df(x,\sbull)\colon E\to F$
is a continuous $\K$-linear map,
for each $x\in U$.
It is known that compositions
of composable $C^n_\K$-maps
are~$C^n_\K$.
Furthermore, continuous linear (or multilinear)
maps are $C^\infty_\K$
(see \cite[Chapter~1]{GaN},
or \cite{RES} for all of this).
\begin{rem}
Keller's $C^n_c$-theory
is used as the basis of infinite-dimensional
Lie theory by many authors
(see \cite{FUN}, \cite{DL3},
\cite{GaN}, \cite{Mil},
\cite{Ne1}, \cite{Wur}).
Others use the
``convenient calculus''~\cite{KaM}.\\[2.5mm]
For some purposes,
it is useful to impose
certain completeness properties
on the locally convex space~$F$
involved.
These are, in decreasing order of strength:
Completeness (every Cauchy net converges);
quasi-completeness (every bounded
Cauchy net converges);
sequential completeness
(every Cauchy sequence converges);
and Mackey completeness
(every Mackey-Cauchy sequence
converges, or equivalently:
the Riemann integral
$\int_0^1\gamma(t)\, dt$ exists in~$F$,
for each smooth curve $\gamma\colon \R\to F$;
see \cite[Theorem~2.14]{KaM} for
further information).
\end{rem}
%
%
\begin{rem}\label{cxanacompa}
If $\K=\C$,
then a map $f\colon E\supseteq U\to F$
is $C^\infty_\C$ if and only if it is
\emph{complex analytic}
in the usual sense (as in~\cite{BaS}),
i.e., $f$ is continuous
and for each $x\in U$, there
exists a $0$-neighbourhood $Y\sub U-x$
and continuous homogeneous
polynomials $p_n\colon E\to F$ of degree~$n$
such that
$f(x+y)=\bigcup_{n=0}^\infty p_n(y)$
for all $y\in Y$.
Such maps are also called holomorphic.
If~$F$ is Mackey complete,
then $f$ is~$C^1_\C$ if and only if it is~$C^\infty_\C$
(see \cite[Propositions~7.4 and~7.7]{BGN}
or \cite[Chapter~1]{GaN} for all of this;
cf.\ \cite{RES}).
For suitable non-Mackey complete~$F$, there
are $C^n_\C$-maps $\C\to F$
for all $n\in \N$ which
are not $C^{n+1}_\C$~(\cite{COM},~\cite{Gro}).
\end{rem}
\begin{rem}
If $f\colon U\to F$ is a map
from an open subset of~$\K$
to a locally convex space,
then $f$ is $C^1_\K$ in the above
sense
if and only if the (real, resp.\ complex)
derivative $f^{(1)}(x)=f'(x)=\frac{df}{dx}(x)$
exists for each $x\in U$, and $f'\colon U\to F$
is continuous.
Likewise, $f$ is $C^n_\K$
if it has continuous
derivatives $f^{(k)}\colon U\to F$
for all $k\in \N_0$ such that $k\leq n$
(where $f^{(k)}:=(f^{(k-1)})'$).
This is easy to see (and spelled out in \cite[Chapter~1]{GaN}).
If $\K=\C$ here, then complex analyticity
of~$f$ simply means that $f$ can be expressed
in the form
$f(z)=\sum_{n=0}^\infty (z-z_0)^na_n$
close to each given point $z_0\in U$,
for suitable elements $a_n\in F$.
\end{rem}
%
%
\begin{rem}\label{complements}
Consider a map
$f\colon U\to F$
from an open set $U\sub \C$ to a Mackey complete
locally convex space~$F$.
Replacing sequential completeness
with Mackey completeness
in \cite[Theorem~3.1]{BaS}
and its proof,
one finds that also
each of the following conditions
is equivalent
to $f$ being a $C^\infty_\C$-map
(see also \cite[Chapter~II]{Gro},
notably Theorems~2.2, 2.3 and~5.5):\footnote{Considering $f$
as a map into the completion
of~$F$, we see that
(a), (b) and (c) remain
equivalent if $F$ is not Mackey complete.
But (a)--(c)
do not imply that $f$ is~$C^1_\C$
(\cite[Example~II.2.3]{Gro},
\cite[Theorem~1.1]{COM}).}
\begin{itemize}
\item[(a)]
$f$ is weakly holomorphic,
i.e.\ $\lambda\circ f\colon U\to \C$
is holomorphic for each $\lambda\in F'$;
\item[(b)]
$\int_{\partial\Delta} f(\zeta)\,d\zeta=0$
for each triangle $\Delta\sub U$;
\item[(c)]
$f(z)=\frac{1}{2\pi i}\int_{|\zeta-z_0|=r}\frac{f(\zeta)}{\zeta-z}\,d\zeta$
for each $z_0\in U$,
$r>0$ such that
$z_0+r\bD\sub U$, and each~$z$
in the interior of the disk $z_0+r\bD$.
\end{itemize}
\end{rem}
\begin{defn}
Given locally convex spaces~$E$ and~$F$,
we let $L(E,F)$ be the vector space
of all continuous linear maps $A\colon E\to F$.
If $\cS$ is a set of bounded
subsets of~$E$, we write
$L(E,F)_\cS$ for $L(E,F)$, equipped with the
topology of uniform convergence
on the sets $M\in \cS$.
Finite intersections of sets
of the form
\[
\lfloor M,U\rfloor\; :=\;
\{A\in L(E,F)\colon A(M)\sub U\}
\]
(for $M\in \cS$ and $U\sub F$ a $0$-neighbourhood)
form a basis for the filter of
$0$-neighbourhoods
of this vector topology.
See \cite[Chapter~III, \S3]{BTV}
for further information.
Given $M\sub E$
and $N\sub E'$, we write
$M^\circ:=\lfloor M,\bD\rfloor\sub E'$
and ${}^\circ N:=\{x\in E\colon
(\forall\lambda\in N)\; \lambda(x)\in \bD\}$
for the polar in~$E'$ (resp., in~$E$).
\end{defn}
\begin{rem}
If~$F$ is Hausdorff
and $F\not=\{0\}$, then $L(E,F)_\cS$
is Hausdorff
if and only $\bigcup_{M\in \cS}M$
is total in~$E$, i.e.,
it spans a dense vector subspace.\\[2.5mm]
In fact, totality of $\bigcup\cS$ is sufficient for
the Hausdorff property
by Proposition~3 in \cite[Chapter~III, \S3, no.\,2]{BTV}.
If $V:=\Spann_\K(\bigcup\cS)$
is not dense in~$E$,
the Hahn-Banach Theorem
provides a linear functional
$0\not= \lambda\in E'$ such that $\lambda|_V=0$.
We pick $0\not = y\in F$. Then $y\tensor \lambda\in
\lfloor M,U\rfloor$
for each $M\in \cS$ and $0$-neighbourhood $U\sub F$,
whence $y\tensor \lambda\in W$
for each $0$-neighbourhood $W\sub L(E,F)_\cS$.
Since $y\tensor\lambda\not=0$,
$L(E,F)_\cS$ is not Hausdorff.
\end{rem}
%
%
\begin{prop}\label{prehypo}
Given a separately
continuous bilinear
map $\beta\colon E_1\times E_2\to F$
and a set $\cS$ of bounded subsets
of~$E_2$,
consider the following conditions:
\begin{itemize}
\item[\rm(a)]
For each $M \in \cS$
and each $0$-neighbourhood $W \sub F$,
there exists a $0$-neighbourhood
$V\sub E_1$
such that $\beta(V\times M)\sub W$.
\item[\rm(b)]
The mapping $\beta^\vee\colon
E_1\to L(E_2,F)_\cS$, $x\mto \beta(x,\sbull)$
is continuous.
\item[\rm(c)]
$\beta|_{E_1\times M}\colon E_1\times M\to F$
is continuous, for each $M\in \cS$.
\end{itemize}
Then {\rm (a)} and {\rm (b)}
are equivalent,
and {\rm (a)} implies {\rm (c)}.
If
%
\begin{equation}\label{simplf}
(\forall M\in \cS)\;(\exists N\in \cS)\;\;\;
\bD M\sub N\,,
\end{equation}
then all of {\rm (a)--(c)}
are equivalent.
\end{prop}
%
%
\begin{defn}\label{defhypo}
A bilinear map $\beta$
which is separately
continuous and satisfies the equivalent conditions
(a) and~(b) of Proposition~\ref{prehypo}
is called \emph{$\cS$-hypocontinuous}
(in the second argument),
or simply \emph{hypocontinuous}
if~$\cS$ is clear from the context.
Hypocontinuity in the first
argument with respect to
a set of bounded subsets of~$E_1$
is defined analogously.
\end{defn}
{\bf Proof of Proposition~\ref{prehypo}.}
For the equivalence
(a)$\aeq$(b)
and the implication (b)$\impl$(c),
see Proposition~3 and~4
in \cite[Chapter~III, \S5, no.\,3]{BTV},
respectively.\\[2.5mm]
We now show that
(c)$\impl$(a) if~(\ref{simplf}) is satisfied.
Given $M\in \cS$ and $0$-neighbourhood
$W\sub F$, by hypothesis
we can find $N\in \cS$ such that
$\bD M\sub N$.
By continuity
of $\beta|_{E_1\times N}$,
there exist
$0$-neighbourhoods~$V$ in~$E_1$ and~$U$ in~$E_2$
such that $\beta(V\times (N\cap U))\sub W$.
Since~$M$ is bounded,
$M \sub nU$ for some $n\in \N$.
Then $\frac{1}{n}M\sub N\cap U$.
Using that~$\beta$ is bilinear, we obtain
$\beta((\frac{1}{n}V)\times M)=\beta(V \times(\frac{1}{n}M))
\sub \beta(V\times (N\cap U))\sub W$.\,\Punkt
%
%
%
\begin{rem}\label{rembalanc}
By Proposition~\ref{prehypo}\,(b),
$\cS$-hypocontinuity
of $\beta\colon E_1\times E_2\to F$
only depends on the topology
on $L(E_2,F)_\cS$, not on~$\cS$ itself.
Given~$\cS$, define
$\cS':=\{\bD M\colon M\in \cS\}$.
Then the topologies
on $L(E_2,F)_\cS$ and $L(E_2,F)_{\cS'}$
coincide (as is clear),
and hence $\beta$
is $\cS$-hypocontinuous
if and only if~$\beta$
is $\cS'$-hypocontinuous.
After replacing~$\cS$ with~$\cS'$,
we can therefore
always assume that~(\ref{simplf})
is satisfied,
whenever this is convenient.
\end{rem}
Each continuous bilinear map is
hypocontinuous
(as condition~(a) in Proposition~\ref{prehypo}
is easy to check), but the converse
is false.
The next proposition compiles
various useful facts.
%
%
\begin{prop}\label{somefcts}
Let $\beta\colon E_1\times E_2\to F$
be an $\cS$-hypocontinuous
bilinear map,
for some set~$\cS$ of bounded
subsets of~$E_2$.
Then the following holds.
\begin{itemize}
\item[\rm(a)]
$\beta(B\times M)$ is bounded
in~$F$, for each bounded subset
$B\sub E_1$ and each $M\in \cS$.
\item[\rm(b)]
Assume that, for each
convergent sequence $(y_n)_{n\in \N}$
in~$E_2$, with limit~$y$,
there exists $M\in \cS$ such that
$\{y_n\colon n\in \N\}\cup\{y\}\sub M$.
Then $\beta$ is sequentially continuous.
\end{itemize}
The condition described
in {\rm (b)} is satisfied,
for example, if $\cS$
is the set of all bounded subsets
of~$E_2$, or the set
of all compact subsets
of~$E_2$.
\end{prop}
\begin{proof}
(a) See Proposition~4
in \cite[Chapter~III, \S5, no.\,3]{BTV}.

(b) See \cite[p.\,157, Remark following \S40,\,1.,\,(5)]{Koe}.
\end{proof}
In many cases,
separately continuous
bilinear maps are automatically
hypocontinuous.
To make this precise,
we recall that
a subset~$B$ of a locally convex space~$E$
is called a \emph{barrel}
if it is closed,
convex, balanced and absorbing.
The space $E$ is called
\emph{barrelled}
if every barrel is a $0$-neighbourhood.
See Proposition~6 in \cite[Chapter~III, \S5, no.\,3]{BTV}
for the following~fact:
%
%
\begin{prop}\label{auothypo}
If $\beta \colon E_1\times E_2\to F$
is a separately continuous
bilinear map and~$E_1$ is barrelled, then~$\beta$
is hypocontinuous with respect to
any set~$\cS$ of bounded subsets of~$E_2$.\,\Punkt
\end{prop}
A simple fact will be useful.
%
%
\begin{la}\label{pla}
Let~$M$ be a topological space, $F$ be a locally convex space,
and $BC(M,F)$ be the space of bounded $F$-valued
continuous functions on~$M$,
equipped with the topology of uniform convergence.
Then the evaluation map
$\mu\colon BC(M,F)\times M\to F$,
$\mu(f,x):=f(x)$
is continuous.
\end{la}
\begin{proof}
Let $(f_\alpha,x_\alpha)$ be a convergent net in
$BC(M,F)\times X$, convergent to $(f,x)$, say.
Then $\mu(f_\alpha,x_\alpha)-\mu(f,x)=
(f_\alpha(x_\alpha)-f(x_\alpha))+(f(x_\alpha)-f(x))$,
where $f_\alpha(x_\alpha)-f(x_\alpha)\to 0$
as $f_\alpha\to f$ uniformly
and $f(x_\alpha)-f(x)\to 0$ as $f$ is continuous.
\end{proof}
We now turn to two paradigmatic
bilinear maps, namely evaluation
and composition.
%
%
\begin{prop}\label{resulteval}
Let~$E$ and~$F$ be locally convex spaces
and~$\cS$ be a set of bounded
subsets of~$E$
which covers~$E$,
i.e., $\bigcup_{M\in \cS}M=E$.
Then the
evaluation map
\[
\ve\colon L(E,F)_\cS \times E\to F\, ,\quad
\ve(A, x)\, :=\, A(x)
\]
is hypocontinuous
in the second argument
with respect to $\cS$.
If~$E$ is barrelled, then~$\ve$ is also hypocontinuous
in the first argument,
with respect to any locally
convex topology~$\cO$
on $L(E,F)$ which is finer
than the topology of pointwise
convergence,
and any set~$\cT$ of bounded subsets
of $(L(E,F),\cO)$.
\end{prop}
\begin{proof}
By Remark~\ref{rembalanc},
we may assume that
$\cS$ satisfies~(\ref{simplf}).
Given $A \in L(E,F)$, we have
$\ve(A,\sbull)=A$, whence~$\ve$
is continuous in the second argument.
It is also continuous in the first argument, as
the topology on $L(E,F)_\cS$
is finer than the topology
of pointwise convergence,
by the hypothesis on~$\cS$.
Let~$M\in \cS$ now.
As $L(E,F)$ is equipped with the
topology of uniform convergence on the sets
in~$\cS$,
the restriction map $\rho\colon L(E,F)\to BC(M,F)$,
$A \mto A|_M$ is continuous.
By Lemma~\ref{pla},
the evaluation map $\mu\colon BC(M,F)\times M\to F$
is continuous. Now $\ve|_{L(E,F)\times M}=\mu\circ (\rho\times \id_M)$
shows that~$\ve|_{L(E,F)\times M}$ is continuous.
Since we assume~(\ref{simplf}),
the implication
``(c)$\impl$(a)'' in
Proposition~\ref{prehypo}
shows that~$\ve$ is $\cS$-hypocontinuous.\\[2.5mm]
Since~$\cO$ is finer than the topology of
pointwise convergence,
the map~$\ve$ remains separately continuous
in the situation described at the end
of the proposition.
Hence, if~$E$ is barrelled,
Proposition~\ref{auothypo}
ensures hypocontinuity
with respect to~$\cT$.
\end{proof}
While it was sufficient
so far to consider an individual
set $\cS$ of bounded subsets
of a given locally convex space,
we now frequently
wish to select such a set $\cS$
simultaneously for each space.
The following definition
captures the situations of interest.
%
%
\begin{defn}\label{bsetfun}
A \emph{bounded set functor}
is a functor $\cS$ from the
category of locally convex spaces
to the category of sets,
with the following properties:
\begin{itemize}
\item[(a)]
$\cS(E)$ is a set of bounded
subsets of~$E$,
for each locally convex space~$E$.
\item[(b)]
If $A\colon E\to F$ is a continuous linear
map, then $A(M)\in \cS(F)$
for each $M\in \cS(E)$,
and $\cS(A)\colon \cS(E)\to\cS(F)$
is the map taking
$M\in \cS(E)$
to its image $A(M)$ under~$A$.
\end{itemize}
Given a bounded set functor~$\cS$
and locally convex spaces~$E$ and~$F$,
we write $L(E,F)_\cS$
as a shorthand for $L(E,F)_{\cS(E)}$.
Also, an $\cS(E_2)$-hypocontinuous
bilinear map $\beta\colon E_1\times E_2\to F$
will simply be called $\cS$-hypocontinuous in the
second argument.
\end{defn}
%
%
\begin{example}\label{threemain}
Bounded
set functors are
obtained if $\cS(E)$
denotes the set of all bounded,
(quasi-)\,compact, or finite
subsets of~$E$, respectively.
We then write $b,c$, resp., $p$ for~$\cS$.
\end{example}
Further examples abound:
For instance, we can let $\cS(E)$
be the set of precompact subsets
of~$E$, or the set of metrizable
compact subsets
(if only Hausdorff spaces are considered).
\begin{rem}
If $\cS$ is a bounded set functor,
and $A\colon E\to F$ a continuous
linear map between locally convex spaces,
then also its adjoint $A'\colon F'_\cS\to E'_\cS$,
$\lambda\mto \lambda\circ A$
is continuous,
because $A'(\lfloor A(M),U\rfloor)\sub \lfloor
M,U\rfloor$ for each $M\in \cS(E)$
and $0$-neighbourhood $U\sub \K$.
\end{rem}
The double use of $f'(x)$
(for differentials) and~$A'$ (for adjoints)
should not cause confusion.
\begin{defn}
If $\cS(E)$ contains
all finite subsets of~$E$,
then $\eta_E(x)\colon E'_\cS\to \K$, $\lambda\mto \lambda(x)$
is continuous for each $x\in E$
and we obtain a linear map
$\eta_E\colon E\to (E'_\cS)'$,
called the \emph{evaluation
homomorphism}.
We say that~$E$ is \emph{$\cS$-reflexive}
if $\eta_E\colon E\to (E'_\cS)'_\cS$
is an isomorphism of topological vector spaces.
If $\cS=b$, we simply speak of a \emph{reflexive}
space; if $\cS=c$, we speak of a \emph{Pontryagin
reflexive} space. Occasionally,
we call $E'_b$ the \emph{strong}
dual of~$E$.
\end{defn}
See Proposition~9 in \cite[Chapter~III, \S5, no.\,5]{BTV}
for the following fact
in the three cases described
in Example~\ref{threemain}.
It might also be deduced from
\cite[\S40,\,5.,\,(6)]{Koe}.
%
%
\begin{prop}\label{resultcomp}
Let $E$, $F$, and $G$
be locally convex spaces
and~$\cS$ be a bounded set functor
such that $\cS(E)$ covers~$E$ and
%
\begin{equation}\label{strnge}
\forall M\in \cS(L(E,F)_\cS)\;\, \forall N \in \cS(E)\;\,
\exists K\in \cS(F)\colon \quad
\ve(M\times N)\sub K\,,
\end{equation}
where $\ve\colon L(E,F)\times E\to F$, $(A,x)\mto A(x)$.
Then the composition map
\[
\Gamma\colon L(F,G)_\cS \times L(E,F)_\cS \to L(E,G)_\cS\,,\quad
\Gamma(\alpha,\beta):=\alpha\circ \beta
\]
is $\cS(L(E,F)_\cS)$-hypocontinuous in the second
argument.\,\Punkt
\end{prop}
\begin{rem}
If $\cS=b$,
then condition~(\ref{strnge})
is satisfied by Proposition~\ref{somefcts}\,(a).
If~$\cS=p$, then $\ve(M\times N)$ is finite
and thus (\ref{strnge}) holds.
If $\cS=c$,
then $\ve|_{M\times N}$
is continuous since
$\ve\colon L(E,F)_\cS\times E\to F$
is $\cS(E)$-hypocontinuous
by Proposition~\ref{resulteval}.
Hence $\ve(M\times N)$
is compact
(and thus (\ref{strnge})
is satisfied).
\end{rem}
{\bf Proof of Proposition~\ref{resultcomp}.}
\emph{$\Gamma$ is continuous in the second argument}:
Let $M\in \cS(E)$, $U\sub G$ be
a $0$-neighbourhood, and $A\in L(F,G)$.
Then $A^{-1}(U)$ is a $0$-neighbourhood in~$F$.
For $B \in L(E,F)$, we have
$A(B(M))\sub U$ if and only if
$B(M)\sub A^{-1}(U)$,
showing that $\Gamma(A,\lfloor M,A^{-1}(U)\rfloor)
\sub \lfloor M,U\rfloor$. Hence, $\Gamma(A,\sbull)$
being linear, it is continuous.\\[2.5mm]
\emph{Continuity in the first argument}:
Let $U\sub G$ be a $0$-neighbourhood,
$B \in L(E,F)$ and $M\in \cS(E)$.
Then $B(M)\in \cS(F)$ by Definition~\ref{bsetfun}\,(b)
and $\Gamma(\lfloor B(M),U\rfloor,B)\sub \lfloor
M,U\rfloor$.\\[2.5mm]
To complete the proof,
let $M\in \cS(L(E,F)_\cS)$,
$U\sub G$ be a $0$-neighbourhood
and~$N\in \cS(E)$.
By~(\ref{strnge}),
there exists $K\in \cS(F)$
such that
$\ve(M\times N)\sub K$.
Note that for all $B\in M$ and
$A\in \lfloor K ,U\rfloor$,
we have $\Gamma(A,B).N=(A\circ B)(N)
=A(B(N))\sub A(K)\sub U$.
Thus
$\Gamma(\lfloor K ,U\rfloor\times M)\sub \lfloor N,U\rfloor$.
Since $\lfloor K,U\rfloor$
is a $0$-neighbourhood in $L(F,G)_\cS$,
condition~(a)
of Proposition~\ref{prehypo} is satisfied.\vspace{.55mm}\,\Punkt
%
%
\begin{rem}\label{nachtrg}
Despite the hypocontinuity of the composition map~$\Gamma$,
it is discontinuous in the
situations specified in the introduction.
To see this, pick
$\lambda\in E'$ and $x\in E$
as described in the introduction.
Let $0\not=z\in G$ and give~$F'$ the topology
induced by $F'\to L(F,G)$, $\zeta\mto z\tensor \zeta$.
There is
$\mu\in G'$ such that $\mu(z)\not=0$.
If~$\Gamma$ was continuous,
then also the following map would be continuous:
$F'\times F\to \K$,
$(\zeta, y)\mto \mu(\Gamma(z\tensor \zeta, y\tensor \lambda)(x))
=\mu(z)\lambda(x)\zeta(y)$.
But this map is a non-zero
multiple of the evaluation map
and hence discontinuous~\cite{Mai}.
\end{rem}
%
%
%
%
%
%
%
%
%
%
%
\section{Differentiability properties
of compositions with\\
hypocontinuous bilinear mappings}\label{mainsec}
In this section, we introduce a new class
of topological spaces (``$k^\infty$-spaces'').
We then discuss compositions
of hypocontinuous bilinear maps
with $C^n$-maps on open subsets
of locally convex spaces
which are $k^\infty$-spaces.\\[2.5mm]
Recall that
a Hausdorff topological space~$X$
is called a \emph{$k$-space}
if, for every subset~$A\sub X$,
the set~$A$ is closed in~$X$
if and only if $A\cap K$ is closed
in~$K$ for each
compact subset~$K\sub X$.
Equivalently, a subset $U\sub X$ is open in~$X$
if and only if
$U\cap K$ is open in~$K$
for each compact subset $K\sub X$.
It is clear that closed subsets,
as well as open subsets
of $k$-spaces
are again $k$-spaces when equipped
with the induced topology.
If $X$ is a $k$-space,
then a map $f\colon X\to Y$ to a topological
space~$Y$ is continuous
if and only if $f|_K\colon K\to Y$
is continuous for each compact subset $K\sub X$,
as is easy to see.\footnote{Given a closed set $A\sub Y$,
the intersection $f^{-1}(A)\cap K=(f|_K)^{-1}(A)$
is closed in~$K$ in the latter case
and thus $f^{-1}(A)$ is closed,
whence $f$ is continuous.}
This property is crucial for the following.
It can also be interpreted
as follows:
$X=\dl_K\, K$\vspace{-.8mm} as a topological space.
%
%
\begin{defn}\label{defkinfty}
We say that a topological space~$X$
is a \emph{$k^\infty$-space}
if it is Hausdorff and
its $n$-fold power
$X^n=X\times\cdots \times X$
is a $k$-space,
for each $n\in \N$.
\end{defn}
%
%
\begin{example}\label{exmetr}
It is well known (and easy to prove)
that every metrizable topological space
is a $k$-space. Finite powers
of metrizable spaces being metrizable,
we see:
\emph{Every metrizable
topological space is a $k^\infty$-space.}
\end{example}
%
%
\begin{example}\label{komeg}
A Hausdorff topological space~$X$
is called a \emph{$k_\omega$-space}
if it is a $k$-space
and \emph{hemicompact},\footnote{An equivalent definition runs as follows:
A Hausdorff space~$X$ is a $k_\omega$-space if and only if
$X=\dl\,K_n$\vspace{-1.6mm}
for an ascending sequence $(K_n)_{n\in \N}$
of compact subsets of~$X$
with union~$X$.}
i.e., there exists a sequence
$K_1\sub K_2\sub \cdots$
of compact subsets
of~$X$ such that $X=\bigcup_{n\in \N}K_n$
and each compact subset of~$X$
is contained in some~$K_n$.
Since finite products
of $k_\omega$-spaces
are $k_\omega$-spaces
(see, e.g., \cite[Proposition~4.2\,(c)]{GGH}),
it follows that each $k_\omega$-space
is a $k^\infty$-space.
For an introduction to $k_\omega$-spaces,
the reader may consult~\cite{GGH}.
\end{example}
%
%
\begin{rem}\label{supplykom}
We remark that $E'_c$ is a $k_\omega$-space
(and hence a $k^\infty$-space),
for each metrizable locally convex space~$E$
(see \cite[Corollary~4.7 and Proposition~5.5]{Aus}).
In particular,
every Silva space~$E$
is a $k_\omega$-space
(and hence a $k^\infty$-space),
i.e., every locally convex direct
limit $E=\dl\,E_n$\vspace{-.7mm}
of an ascending sequence
$E_1\sub E_2\sub\cdots$
of Banach spaces, such that the inclusion
maps $E_n\to E_{n+1}$ are compact
operators
(see \cite[Example~9.4]{DL3}).
\end{rem}
Having set up the terminology,
let us record a simple, but useful
observation.
%
%
\begin{thm}\label{mainobs}
Let $n\in \N_0\cup\{\infty\}$.
If $n=0$, let $U=X$ be a topological
space. If $n\geq 1$,
let $X$ be a locally convex space
and $U\sub X$ be an open subset.
Let
$\beta\colon E_1\times E_2\to F$
be a bilinear map
and $f\colon U\to E_1\times E_2$ be a $C^n$-map.
Assume that at least one of {\rm (a), (b)} holds:
\begin{itemize}
\item[\rm(a)]
$X$ is metrizable
and $\beta$ is sequentially continuous; or:
\item[\rm(b)]
$X$ is a $k^\infty$-space
and $\beta$ is hypocontinuous
in the second argument
with respect to a set~$\cS$ of bounded
subsets of~$E_2$
which contains all compact
subsets of~$E_2$.
\end{itemize}
Then $\beta\circ f\colon U\to F$
is~$C^n$.
\end{thm}
\begin{proof}
It suffices to consider the case where
$n<\infty$.
The proof is by induction.\\[2.5mm]
We assume~(a) first.
If $n=0$, let $(x_k)_{k\in \N}$
be a convergent sequence in~$U$,
with limit~$x$. Then $f(x_k)\to x$
by continuity of~$f$ and hence
$\beta(f(x_k))\to \beta(f(x))$,
since~$\beta$ is sequentially continuous.\\[2.5mm]
Now let $n\geq 1$ and assume that the assertion holds
if~$n$ is replaced with $n-1$.
Given $x\in U$ and $y\in X$,
let $(t_k)_{k\in \N}$
be a sequence in $\K\setminus \{0\}$
such that $x+t_ky\in U$ for
each $k\in \N$, and $\lim_{k\to\infty}\,t_k=0$.
Write $f=(f_1,f_2)$ with $f_j\colon U\to E_j$.
Then
\begin{eqnarray*}
\hspace*{-2mm}\lefteqn{\frac{\beta(f(x+t_ky))-\beta(f(x))}{t_k}}\qquad\\
&=&
\beta\left(\frac{f_1(x+t_ky)-f_1(x)}{t_k},f_2(x+t_ky)\right)
+
\beta\left(f_1(x),\frac{f_2(x+t_ky)-f_2(x)}{t_k}\right)\\
&\to & \beta(df_1(x,y),f_2(x))+\beta(f_1(x),df_2(x,y))\quad
\mbox{as $k\to\infty$,}
\end{eqnarray*}
by continuity of~$f$
and sequential continuity of~$\beta$.
Hence the limit $d(\beta\circ f)(x,y)=
{\displaystyle \lim_{t\to 0}}\,
\frac{\beta(f(x+ty))-\beta(f(x))}{t}$\vspace{-1mm}
exists, and is given by
%
\begin{equation}\label{prodrul}
d(\beta\circ f)(x,y)\;=\;
\beta(df_1(x,y),f_2(x))+\beta(f_1(x),df_2(x,y))\,.
\end{equation}
The mappings
$g_1, g_2 \colon U\times X \to E_1\times E_2$
defined via
$g_1(x,y):=(df_1(x,y),f_2(x))$ and
$g_2(x,y):=(f_1(x),df_2(x,y))$
are $C^{n-1}_\K$.
Since
%
\begin{equation}\label{abstrprod}
d(\beta\circ f)\; =\; \beta\circ g_1+\beta\circ g_2
\end{equation}
by~(\ref{prodrul}),
we deduce from the inductive hypotheses
that $d(\beta\circ f)$ is~$C^{n-1}_\K$ and hence
continuous.
Thus~$\beta\circ f$ is~$C^1_\K$
with $d(\beta\circ f)$ a $C^{n-1}_\K$-map
and hence $\beta\circ f$ is~$C^n_\K$,
which completes the inductive proof in the situation of~(a).\\[2.5mm]
In the situation of~(b),
let $K\sub U$ be compact.
Then $f_2(K)\sub E_2$ is compact
and hence $f_2(K)\in\cS$, by hypothesis.
Since $\beta|_{E_1\times f_2(K)}$
is continuous by Proposition~\ref{prehypo}\,(c),
we see that $(\beta\circ f)|_K=\beta|_{E_1\times f_2(K)}\circ f|_K$
is continuous. Since~$X$ and hence also its open subset~$U$
is a $k$-space, it follows that $\beta\circ f$
is continuous, settling the case $n=0$.\\[2.5mm]
Now let $n\geq 1$ and assume that the assertion holds
if~$n$ is replaced with~$n-1$.
Since
$\beta$ is sequentially continuous
by Proposition~\ref{somefcts}\,(b),
we see as in case~(a) that
the directional derivative
$d(\beta\circ f)(x,y)$ exists,
for all $(x,y)\in U\times X$,
and that $d(\beta\circ f)$ is given by~(\ref{abstrprod}).
Since $g_1$ and~$g_2$ are $C^{n-1}_\K$,
the inductive hypothesis
can be applied to the summands
in~(\ref{abstrprod}).
Thus $d(\beta\circ f)$ is~$C^{n-1}_\K$,
whence $\beta\circ f$ is~$C^1_\K$ with
$d(\beta\circ f)$ a $C^{n-1}_\K$-map,
and so $\beta\circ f$ is~$C^n_\K$.
\end{proof}
Combining Proposition~\ref{resultcomp}
and Theorem~\ref{mainobs},
as a first application
we obtain an affirmative answer
to an open question formulated
by Serge Lang~\cite[p.\,8, Remark]{Lan}.
%
%
\begin{cor}\label{corlang}
Let $U$ be an
open subset of a Fr\'{e}chet space,
$E$, $F$ and~$G$ be Fr\'{e}chet
spaces, and $f\colon U\to L(E,F)_b$, $g\colon U \to L(F,G)_b$
be continuous maps.
Then also the mapping\linebreak
$U\to L(E, G)_b$, $x\mto g(x)\circ f(x)$
is continuous.\Punkt
\end{cor}
%
%
%
%
%
%
%
%
%
%
\section{Holomorphic families
of operators}\label{secholfam}
In this section, we compile
conclusions from the previous results
and some useful additional material.
Specializing to the case $X:=\K:=\C$ and
$n:=\infty$,
we obtain results concerning
holomorphic families
of operators, i.e.,
holomorphic maps $U\to L(E,F)_\cS$,
where $U$ is an open subset
of~$\C$.
Among other things,
such holomorphic families
are of interest
for representation theory and cohomology
(see \cite{BaO} and \cite{BO2}).
%
%
\begin{prop}\label{cor2}
Let
$X$, $E$, $F$ and $G$ be locally convex spaces
over~$\K$, such that $X$ is a $k^\infty$-space
$($e.g.\ $X=\K=\C)$.
Let $U\sub X$ be open, $n\in \N_0\cup\{\infty\}$
and
$f\colon U\to L(E,F)_\cS$ as well
as $g\colon U\to L(F,G)_\cS$
be $C^n_\K$-maps, where 
$\cS=b$ or $\cS=c$.
Then also the map
$U\to L(E, G)_\cS$,
$z\mto g(z)\circ f(z)$
is~$C^n_\K$.
\end{prop}
\begin{proof}
The composition map $L(F,G)_\cS\times L(E,F)_\cS\to L(E,G)_\cS$
is hypocontinuous with
respect to $\cS(L(E,F)_\cS)$,
by Proposition~\ref{resultcomp}.
Hence Theorem~\ref{mainobs}\,(b) applies.
\end{proof}
The remainder of this section
is devoted to the proof of the following result.
Here $\cS$ is a bounded set functor
such that $\cS(E)$
covers~$E$,
for each locally convex space~$E$.
\begin{prop}\label{cor1}
Let $E$, $F$ and~$X$ be locally convex spaces
over~$\K$.
If $\eta_F\colon F\to (F_\cS')'_\cS$
is continuous,
then $g\colon U\to L(F'_\cS,E'_\cS)_\cS$, $z\mto f(z)'$
is $C^n_\K$,
for each $n\in \N_0\cup\{\infty\}$
and $C^n_\K$-map
$f\colon U\to L(E,F)_\cS$ on an open
subset~$U\sub X$.
\end{prop}
The proof of
Proposition~\ref{cor1} exploits
the continuity of the formation of
adjoints.
\begin{prop}\label{adj}
Let $E$, $F$ be locally convex 
spaces and $\cS$ be a bounded set functor
such that $\cS(F)$ covers~$F$.
If the evaluation
homomorphism $\eta_F\colon F\to (F'_\cS)'_\cS$ is continuous,
then
%
\begin{equation}\label{dfnPsi}
\Psi\colon L(E,F)_\cS\to L(F'_\cS,E'_\cS)_\cS\,,
\quad
\alpha\mto \alpha'
\end{equation}
is a continuous linear map.
\end{prop}
\begin{proof}
After replacing $\cS(V)$ with
$\{r_1M_1\cup\cdots\cup r_nM_n\colon r_1,\ldots,r_n\in \K,
M_1,\ldots,M_n\in \cS(V)\}$
for each locally convex space~$V$
(which does not change $\cS$-topologies),
we may assume that
$\cS(E)$ is closed under finite unions
and multiplication with scalars.
Let $M\in \cS(F'_\cS)$ and $U\sub E'_\cS$ be a $0$-neighbourhood;
we have to show that $\Psi^{-1}(\lfloor M,U\rfloor)$ is a
$0$-neighbourhood in $L(E,F)_\cS$.
After shrinking~$U$, without loss of generality
$U=N^\circ$ for some $N\in \cS(E)$
(by our special hypothesis concerning~$\cS(E)$).
For $\alpha\in L(E,F)$, we have
\begin{eqnarray*}
\alpha'\in\lfloor M,U\rfloor & \aeq & (\forall \lambda\in M)\;\;
\lambda\circ\alpha=\alpha'(\lambda)\in U=N^\circ \\
&\aeq & (\forall \lambda\in M)\,(\forall x\in N)\;\;
|\lambda (\alpha(x))|\leq 1\\
&\aeq & \alpha(N)\sub {}^\circ M\\
&\aeq & \alpha\in\lfloor N, {}^\circ M\rfloor.
\end{eqnarray*}
Since $M\in \cS(F'_\cS)$
and $\eta_F$ is continuous,
the polar ${}^\circ M=\eta_F^{-1}(M^\circ)$
is a $0$-neighbourhood in~$F$.
Thus $\lfloor N,{}^\circ M \rfloor =\Psi^{-1}(\lfloor M,U\rfloor)$
is a $0$-neighbourhood in $L(E,F)_\cS$.
\end{proof}
{\bf Proof of Proposition~\ref{cor1}.}
Since $f$ is a $C^n_\K$-map
and $\Psi$ in Proposition~\ref{adj}
is continuous linear and hence a $C^\infty_\K$-map,
also $g=\Psi\circ f$ is~$C^n_\K$.\vspace{2.5mm}\,\Punkt

\noindent
The locally convex spaces~$E$
such that $\eta_E\colon E\to (E_b')'_b$
is continuous are known as ``quasi-barrelled''
spaces.
They can characterized easily.
Recall that a
subset~$A$ of a locally convex space~$E$ is called
\emph{bornivorous} if it absorbs
all bounded subsets of~$E$.
The space~$E$ is called \emph{bornological}
if every convex, balanced, bornivorous
subset of~$E$ is a $0$-neighbourhood.
See Proposition~2 in \cite[\S11.2]{Jar}
(and the lines following it)
for the following simple fact:\\[2.5mm]
\emph{The evaluation homomorphism
$\eta_E\colon E\to (E'_b)'_b$
is continuous
if and only if each closed,
convex, balanced subset~$A\sub E$ which absorbs
all bounded subsets of~$E$ $($i.e., each
bornivorous barrel~$A)$ is a $0$-neighbourhood
in~$E$.}\\[2.5mm]
Thus
$\eta_E\colon E\to (E'_b)'_b$ is continuous if~$E$
is bornological or barrelled.
It is also known
that $\eta_E\colon E\to (E'_c)'_c$
is continuous if~$E$ is a $k$-space
(cf.\ \cite[Corollary~5.12 and Proposition~5.5]{Aus}).\\[2.5mm]
We mention that
Proposition~\ref{cor2} generalizes
\cite[Lemma~2.4]{BaO}, where $X=\K=\C$, $\cS=b$,
$E$ is assumed to be a Montel space and
$E$, $F$ $G$ are complete
(see last line of \cite[p.\,637]{BaO}).
The method of proof used in \emph{loc.\,cit.}\
depends on completeness properties of $L(E,G)_b$,
because the
characterization of holomorphic functions
via Cauchy integrals (as in Remark~\ref{complements}\,(c))
requires Mackey completeness.
Proposition~\ref{cor1}
generalizes \cite[Lemma~2.3]{BaO},
where $X=\K=\C$, $\cS=b$ and continuity of $\eta_F$ is presumed as well
(penultimate sentence of their proof), and
whose proof is valid whenever $L(F',E')_b$
is at least Mackey complete
(since only weak holomorphicity of~$g$ is checked there).
%
%
%
%
%
%
%
%
%
%
%
%
%
%
\section{Locally convex Poisson
vector spaces}\label{secpoisson}
We now consider
locally convex Poisson vector spaces
in a framework which arose from~\cite{GLN}.
Fundamental facts
concerning such spaces
will be proved,
using hypocontinuity as a tool.
%
%
\begin{numba}\label{thesett}
Throughout this
section, we let
$\cS$ be a bounded set functor
such that the following holds
for each locally convex space~$E$:
\begin{itemize}
\item[(a)]
$\cS(E)$ contains
all compact subsets of~$E$; and:
\item[(b)]
For each $M\in \cS(E'_\cS)$
and $N\in \cS(E)$,
the set $\ve(M\times N)\sub \K$
is bounded,
where $\ve\colon E'\times E\to\K$
is the evaluation map.
\end{itemize}
Since $\cS(\K)$ contains all
compact sets and each bounded
subset of~$\K$ is contained
in a compact set,
condition~(b) means
that there exists $K\in \cS(\K)$
such that
$\ve(M\times N)\sub K$.
\end{numba}
%
%
\begin{defn}\label{defnposp}
An \emph{$\cS$-reflexive locally convex Poisson vector space}
is a locally convex space~$E$
which is $\cS$-reflexive
and a $k^\infty$-space,
together with an $\cS$-hypocontinuous
bilinear map
$[.,.] \colon E'_\cS \times E'_\cS
\to E'_\cS$,
$(\lambda,\eta)\mto [\lambda,\eta]$
which makes $E'_\cS$ a Lie algebra.
\end{defn}
Of course, we are mostly interested
in the case where~$[.,.]$ is continuous,
but only $\cS$-hypocontinuity
is required for the basic results
described below.
%
\begin{rem}\label{remtypi}
We mainly have two choices of~$\cS$
in mind.
\begin{itemize}
\item[(a)]
The case~$\cS=b$.
If~$E$ is a Hilbert
space, a reflexive Banach space,
a nuclear
Fr\'{e}chet space, or the strong dual
of a nuclear Fr\'{e}chet space,
then both reflexivity is satisfied
and also the $k^\infty$-property
(by Example~\ref{exmetr}
and Remark~\ref{supplykom}).\footnote{Let~$E$
be a nuclear Fr\'{e}chet space.
Then~$E$ is Pontryagin reflexive
\cite[Propositions~15.2 and 2.3]{Ban}.\linebreak
By \cite[Theorem~16.1]{Ban}
and \cite[Proposition~5.9]{Aus},
also $E'_c$ is nuclear
and Pontryagin reflexive.
Since $E$ is metrizable
and hence a $k$-space, $E'_c$
is complete~\cite[Proposition~4.11]{Aus}.
We now see with \cite[Proposition~50.2]{Tre}
that every closed, bounded subset
of $E$ (and~$E'_c$) is compact.
Hence $E_b'=E_c'$, $(E_b')'_b=(E'_c)'_c$
and both $E$ and $E'_b=E'_c$
are also reflexive.
By Remark~\ref{supplykom}, $E'_b=E'_c$
is a $k^\infty$-space.}
\item[(b)]
If $\cS=c$, then the scope widens
considerably.
For example, every Fr\'{e}chet
space~$E$ is both Pontryagin reflexive
(see \cite[Propositions~15.2 and~2.3]{Ban})
and a $k^\infty$-space
(see Example~\ref{exmetr}),
and the same holds for~$E'_c$
(see \cite[Proposition~5.9]{Aus}
and Remark~\ref{supplykom}).
\end{itemize}
\end{rem}
%
%
%
\begin{rem}\label{remtypl}
The most
typical examples of $\cS$-reflexive
locally convex Poisson vector
spaces are dual spaces of topological
Lie algebras.
More precisely,
let $\cS=b$ or $\cS=c$,
and $(\cg,[.,.]_\cg)$
be a locally convex topological
Lie algebra.
If $\cg$ is $\cS$-reflexive
and $\cg'_\cS$
happens to be a $k^\infty$-space,
then $E:=\cg'_\cS$
is an $\cS$-reflexive
locally convex Poisson vector space
with Lie bracket
defined via
$[\lambda,\mu]:=[\eta_\cg^{-1}(\lambda),
\eta_\cg^{-1}(\mu)]_\cg$
for $\lambda, \mu\in E'=(\cg'_\cS)'$,
using the isomorphism
$\eta_\cg\colon \cg\to (\cg'_\cS)_\cS'$.
Here are typical examples.
\begin{itemize}
\item[(a)]
If~$\cg$ is a Banach-Lie algebra
whose underlying Banach space is
reflexive,
then $\cg'_b$ is a \emph{reflexive} locally convex
Poisson vector space
(i.e., w.r.t.\ $\cS=b$);
see Remark~\ref{remtypi}\,(a).
\item[(b)]
If $\cg$ is a Fr\'{e}chet-Lie algebra
(a topological Lie algebra
which is a Fr\'{e}chet space),
then $\cg'_c$ is
a \emph{Pontryagin
reflexive} locally convex
Poisson vector space (i.e., with respect to $\cS=c$),
by Remark~\ref{remtypi}\,(b).
\item[(c)]
If $\cg$ is a Silva-Lie algebra,
then $\cg$ is reflexive
(hence also Pontryagin reflexive),
and $\cg_b'=\cg'_c$ is
a Fr\'{e}chet-Schwartz space
(see \cite{Flo})
and hence a $k^\infty$-space.
Therefore $\cg'_b=\cg'_c$ is a reflexive and Pontryagin reflexive
locally convex Poisson vector space.
\end{itemize}
If a topological group~$G$
is a projective limit
$\pl\,G_n$\vspace{-.7mm}
of a projective sequence \mbox{$\cdots \to G_2\to G_1$}
of finite-dimensional Lie groups,
then $\cg:=\pl\,L(G_n)\isom \R^\N$\vspace{-.7mm}
can be considered as the Lie algebra
of~$G$ and coadjoint orbits of~$G$ in~$\cg'$
can be studied~\cite{Ne2},
where $\cg'_c$ ($=\cg'_b$)
is a Pontryagin reflexive (and reflexive)
locally convex Poisson vector space,
by~(b).\\[2.5mm]
If a group~$G$
is a union $\bigcup_{n\in \N}\, G_n$
of finite-dimensional
Lie groups $G_1\sub G_2\sub\cdots$,
then~$G$ can be made an infinite-dimensional
Lie group with Lie algebra $\cg=\dl\, L(G_n)\isom \R^{(\N)}$\vspace{-.7mm}
(see \cite{FUN}),
where~$\cg'_c$ ($=\cg'_b$)
is a Pontryagin reflexive (and reflexive)
locally convex Poisson vector space,
by~(c). Again coadjoint orbits can be studied~\cite{GLN}.
Manifold structures on them do not pose
problems, since all homogeneous spaces
of~$G$ are manifolds~\cite[Proposition~7.5]{FUN}.
\end{rem}
Given a Lie algebra $(\cg,[.,.])$
and $x\in \cg$, we write $\ad_x:=\ad(x):=[x,.]\colon
\cg\to\cg$, $y\mto [x,y]$.\\[2.5mm]
Definition~\ref{dfnpoivec}
can be adapted to spaces which are not $\cS$-reflexive,
along the lines of~\cite{OR1},~\cite{OR2}:
\begin{defn}\label{nonrcase}
A \emph{locally convex Poisson vector space}
with respect to~$\cS$
is a locally convex space~$E$
that is a $k^\infty$-space
and whose evaluation homomorphism
$\eta_E\colon E\to (E'_\cS)'_\cS$
is~a topological embedding,
together with an $\cS$-hypocontinuous
bilinear map
\mbox{$[.,.] \colon E'_\cS \times E'_\cS \to E'_\cS$,}
$(\lambda,\eta)\mto [\lambda,\eta]$
which makes $E'_\cS$ a Lie algebra,
and such that
\begin{equation}\label{compcond}
\eta_E(x)\circ \ad_\lambda\; \in \; \eta_E(E)\quad
\mbox{for all $x\in E$ and $\lambda\in E'$.}
\end{equation}
\end{defn}
Identifying~$E$ with $\eta_E(E)\sub (E'_\cS)_\cS'$,
we can rewrite (\ref{compcond}) as
\begin{equation}\label{compcond2}
(\ad_\lambda)'(E)\; \sub\; E\quad
\mbox{for all $\lambda\in E'$.}
\end{equation}
\begin{rem}
Every $\cS$-reflexive
Poisson vector space $(E,[.,.])$
in the sense of Definition~\ref{defnposp}
also is a Poisson vector space with respect to~$\cS$,
in the sense of Definition~\ref{nonrcase}.
In fact,
since $[.,.]$ is separately
continuous, the linear map
$\ad_\lambda=[\lambda ,.]\colon E'_\cS\to E'_\cS$
is continuous, for each $\lambda \in E'_\cS$.
Hence $\alpha\circ \ad_\lambda \in (E_\cS')'=\eta_E(E)$
for each $\alpha \in (E_\cS')'$,
and thus (\ref{compcond}) is satisfied.
\end{rem}
\begin{rem}
If $\cS=b$, then $\eta_E\colon E\to (E'_b)'_b$
is a topological embedding if and only $\eta_E$
is continuous, i.e., if and only if~$E$
is quasi-barrelled in the sense recalled
in Section~\ref{secholfam} (see \cite[\S11.2]{Jar}).
Most locally convex spaces of practical interest
are bornological or barrelled
and hence quasi-barrelled.\\[2.5mm]
If $\cS=c$, then $\eta_E$ is a topological embedding
\emph{automatically} in the situation of Definition~\ref{nonrcase}
as we assume that $E$ is a $k^\infty$-space
(and hence a $k$-space).
In fact, $\eta_E\colon E\to (E'_c)'_c$
is injective (by the Hahn-Banach theorem)
for each locally convex space~$E$,
and open onto its image
(cf.\ \cite[Proposition~6.10]{Aus}
or \cite[Lemma~14.3]{Ban}).
Hence $\eta_E\colon E\to (E'_c)'_c$
is an embedding if and only
if it is continuous,
which holds if~$E$ is a $k$-space
(cf.\ \cite[Lemma~14.4]{Ban}).
\end{rem}
\begin{rem}
Since reflexive Banach spaces
are rather rare,
the more complicated
non-reflexive theory
cannot be avoided
in the study of Banach-Lie-Poisson
vector spaces
(as in~\cite{OR1}, \cite{OR2}).
By contrast,
typical non-Banach locally convex spaces
are reflexive
and hence fall within the simple,
basic framework of Definition~\ref{defnposp}.
And the class of Pontryagin reflexive spaces
is even more comprehensive.
\end{rem}
%
%
\begin{defn}\label{dfnpoivec}
Let $(E,[.,.])$ be a locally
convex Poisson vector space
with respect to~$\cS$,
and $U\sub E$ be open.
Given
$f,g\in C^\infty_\K(U,\K)$,
we define a function
$\{f,g\}\colon U\to\K$ via
%
\begin{equation}\label{poissonbr}
\{f,g\}(x)\; :=\; \langle [f'(x),g'(x)], x\rangle\quad
\mbox{for $x\in U$,}
\end{equation}
where $\langle .,.\rangle\colon E'\times E\to \K$,
$\langle \lambda,x\rangle :=\lambda(x)$
is the evaluation map and $f'(x)=df(x,.)$.\\[2.5mm]
Condition (\ref{compcond}) in Definition~\ref{nonrcase}
enables us to define
a map $X_f\colon U \to E$ via
\begin{equation}\label{defhamil}
X_f(x)\, :=\,  \eta_E^{-1}\bigl(\eta_E(x)
\circ \ad(f'(x))\bigr)\quad
\mbox{for $\,x\in U$,}
\end{equation}
where $\eta_E\colon E\to (E'_\cS)'_\cS$
is the evaluation homomorphism.
\end{defn}
%
%
\begin{thm}\label{mainappoi}
Let $(E,[.,.])$
be a locally
convex Poisson vector space
with respect to~$\cS$
and $U\sub E$ be an open subset.
Then
\begin{itemize}
\item[\rm(a)]
$\{f,g\}\in C^\infty_\K(U,\K)$,
for all $f,g\in C^\infty_\K(U,\K)$.
\item[\rm(b)]
For each $f\in C^\infty_\K(U,\K)$, the
map $X_f\colon U \to E$ is $C^\infty_\K$.
\end{itemize}
\end{thm}
The following fact will help us to prove
Theorem~\ref{mainappoi}.
%
%
\begin{la}\label{frombdls}
Let $E$ and~$F$ be locally convex spaces,
$U\sub E$ be open and $f\colon U\to F$ be
a $C^\infty_\K$-map.
Then also the map
$f'\colon U\to L(E,F)_\cS$,
$x\mto f'(x)=df(x,\sbull)$
is $C^\infty_\K$, for each\linebreak
set $\cS$ of bounded subsets
of~$E$ such that $L(E,F)_\cS$
is Hausdorff.
\end{la}
\begin{proof}
For $\cS=b$, see~\cite{HYP}.
The general case is a trivial consequence
of the case~$\cS=b$.
\end{proof}
{\bf Proof of Theorem~\ref{mainappoi}.}
(a) The maps $f'\colon U\to L(E,\K)_\cS=E'_\cS$
and $g'\colon U\to E'_\cS$ are $C^\infty_\K$
by Lemma~\ref{frombdls},
$E$ is a $k^\infty$-space by hypothesis,
and~$[.,.]$ is $\cS$-hypocontinuous.
Hence
$h:= [.,.]\circ (f',g')\colon U\to E'_\cS$,
$x\mto [f'(x),g'(x)]$ is $C^\infty_\K$,
by Theorem~\ref{mainobs}\,(b).
The evaluation map
$\ve\colon E'_\cS\times E\to \K$
is $\cS$-hypocontinuous
in the second argument
by Proposition~\ref{resulteval},
and the inclusion map $\iota\colon U\to E$ is~$C^\infty_\K$.
Hence $\{f,g\}=\ve\circ (h,\iota)$
is $C^\infty_\K$, by Theorem~\ref{mainobs}\,(b).\vspace{1mm}

(b) Since the bilinear map
$[.,.]\colon E'_\cS\times E'_\cS\to E'_\cS$
is $\cS$-hypocontinuous, the linear map
$[.,.]^\vee\colon E'_\cS\to L(E'_\cS,E'_\cS)_\cS$,
$\lambda\mto [\lambda ,.]=\ad(\lambda)$
is continuous, by Proposition~\ref{prehypo}\,(b).
Hence $h\colon U \to L(E'_\cS,E'_\cS)_\cS$, $h(x):=
\ad(f'(x))$ is $C^\infty_\K$,
using Lemma~\ref{frombdls}.
The composition map
$\Gamma\colon (E'_\cS)'_\cS\times L(E'_\cS,E'_\cS)_\cS\to
(E'_\cS)'_\cS$, $(\alpha ,A)\mto \alpha\circ A$
is $\cS$-hypocontinuous in the second argument,
because condition~(b) in \S\,\ref{thesett}
ensures that
Proposition~\ref{resultcomp}
can be applied.
Then $V:=\{A\in L(E'_\cS,E'_\cS)_\cS\colon \mbox{$(\forall
x\in E)\; \eta_E(x)\circ A\in \eta_E(E)$}\}$
is a vector subspace of $L(E'_\cS,E'_\cS)_\cS$
and the bilinear map $\Theta\colon E\times V\to E$,
$\Theta(x,A):=\eta_E^{-1}(\Gamma(\eta_E(x),A))$
is $\cS$-hypocontinuous
in its second argument,
using that $\eta_E$ is
an isomorphism of topological vector spaces
onto its image. The inclusion map
$\iota\colon U\to E$, $x\mto x$ being
$C^\infty_\K$,
Theorem~\ref{mainobs} shows that
$X_f=\Theta \circ (\iota ,h)$
is~$C^\infty_\K$.\,\Punkt
%
%
%
\begin{defn}\label{defhampost}
Let $(E,[.,.])$ be a locally convex
Poisson vector space with respect to~$\cS$,
and $U\sub E$ be an open subset.
\begin{itemize}
\item[(a)]
The map $\{.,.\}\colon C^\infty_\K(U,\K)\times
C^\infty_\K(U,\K)\to C^\infty_\K(U,\K)$
taking $(f,g)$ to $\{f,g\}$
(as in Definition~\ref{dfnpoivec})
is called the \emph{Poisson bracket}.
\item[(b)]
Given $f\in C^\infty_\K(U,\K)$,
the map $X_f\colon U\to E$
is a smooth vector field on~$U$,
by Theorem~\ref{mainappoi}\,(b).
It is called the \emph{Hamiltonian vector field}
associated with~$f$.
\end{itemize}
\end{defn}
%
%
\begin{rem}\label{ispbr}
Basic differentiation rules entail
that $\{.,.\}$ makes
$C^\infty_\K(U,\K)$
a Poisson algebra,
i.e.,
$(C^\infty_\K(U,\K),\{.,.\})$
is a Lie algebra\footnote{The Jacobi identity can be established
as in the proof of \cite[Theorem~4.2]{OR1}.}
and $\{f,.\}\colon$
$C^\infty_\K(U,\K)\to
C^\infty_\K(U,\K)$, $g\mto \{f,g\}$
is a derivation for the
commutative,
associative $\K$-algebra
$(C^\infty_\K(U,\K),\cdot)$,
for each $f\in C^\infty_\K(U,\K)$.
\end{rem}
%
%
%
%
%
%
%
%
%
%
%
%
%
%
\section{Continuity properties of the Poisson bracket}\label{nwsec1}
If $E$ and~$F$ are locally convex spaces,
$U\sub E$ is an open set
and $n\in \N_0\cup\{\infty\}$,
then $C^n_\K(U,F)$
carries a natural topology
(the ``$C^n$-topology''),
namely the initial
topology with respect to the maps
\[
C^\infty_\K(U,F)\to C(U\times E^k,F)_{c.o.}\,\quad
f \mto d^kf
\]
for $k\in \N_0$ such that $k\leq n$,
where the right hand side is
equipped with the compact-open topology,
$d^0f:=f$
and $d^kf(x,y_1,\ldots, y_k):=
(D_{y_k}\cdots D_{y_1}f)(x)$
is defined as an iterated
directional derivative, if $k\geq 1$.
Our goal is the following result:
%
%
\begin{thm}\label{poiss2}
Let $(E,[.,.])$ be a locally
convex Poisson vector space with respect to
$\cS=c$.
Let $U\sub E$ be open.
Then the Poisson bracket
\[
\{.,.\}\colon C^\infty_\K(U,\K)\times
C^\infty_\K(U,\K)\to C^\infty_\K(U,\K)
\]
is hypocontinuous with respect to compact subsets
of $C^\infty_\K(U,\K)$.
If $[.,.]\colon E'_c\times E'_c\to E'_c$
is continuous, then also
the Poisson bracket
is continuous.
\end{thm}
\begin{rem}
The topology on spaces
of smooth maps
goes along well with the $\cS$-topology
on spaces of operators if $\cS=c$,
since it enables to control
smooth functions and their differentials
on compact sets.
For $\cS\not=c$ (notably, for $\cS=b$),
there is no clear connection between
the topologies,
and one cannot hope
for an analogue of Theorem~\ref{poiss2}.
\end{rem}
Various auxiliary results are
needed to prove
Theorem~\ref{poiss2}.
With little risk of confusion
with subsets of spaces
of operators,
given a $0$-neighbourhood $W\sub F$
and a compact set $K\sub U$ we shall write
$\lfloor K, W\rfloor:=\{f \in C(U,F)\colon
f(K)\sub W\}$.
%
%
\begin{la}\label{hlp0}
Let $E, F$ be locally convex spaces
and $U\sub E$ be open.
Then the linear map
\[
D\colon C^\infty_\K(U,F)\to C^\infty_\K(U, L(E,F)_c)\,,\quad
f\mto f'
\]
is continuous.
\end{la}
\begin{proof}
The map $D$ is linear
and also $C^\infty(U,L(E,F)_c)\to
C(U\times E^k,L(E,F)_c)$,
$f \mto d^k f$ is linear,
for each $k\in \N_0$.
Hence
\begin{equation}\label{djD}
d^k\circ D\colon C^\infty(U,F)\to C(U\times E^k,L(E,F)_c)_{c.o.}
\end{equation}
is linear,
whence it will be continuous if it is continuous
at~$0$.
We pick a typical $0$-neighbourhood
in $C(U\times E^k,L(E,F)_c)_{c.o.}$,
say $\lfloor K,V\rfloor$
with a compact subset $K\sub U\times E^k$
and a $0$-neighbourhood $V\sub L(E,F)_c$.
After shrinking~$V$, we may assume
that $V=\lfloor A ,W\rfloor$
for some compact set $A\sub E$ and $0$-neighbourhood
$W\sub F$.\\[2.5mm]
We now recall that for $f\in C^\infty_\K(U,F)$,
we have
%
\begin{equation}\label{usfrec}
d^k(f')(x,y_1,\ldots,y_k)
\; =\; d^{k+1}f(x,y_1,\ldots, y_k,\sbull)\colon E\to F
\end{equation}
for all $k\in \N_0$, $x\in U$ and
$y_1,\ldots, y_k\in E$
(see~\cite{HYP}).
Since $\lfloor K\times A, W\rfloor$
is an open $0$-neighbourhood in $C(U\times E^{k+1},F)$
and the map $C^\infty(U,F)
\to C(U\times E^{k+1},F)_{c.o.}$, $f \mto d^{k+1}f$
is continuous, we see that the set $\Omega$
of all $f \in C^\infty(U,F)$
such that $d^{k+1}f \in\lfloor K\times A, W\rfloor$
is a $0$-neighbourhood in
$C^\infty(U,F)$.
In view of (\ref{usfrec}),
we have $d^k(f')\in \lfloor K,\lfloor A,W\rfloor\rfloor$
for each $f \in\Omega$.
Hence $d^k\circ D$ from~(\ref{djD})
is continuous at~$0$, as required.
\end{proof}
%
%
\begin{la}\label{wannerf}
Let $X$ be a Hausdorff topological space,
$F$ be a Hausdorff locally convex space,
$K\sub X$ be compact
and also $M\sub C(X,F)_{c.o.}$ be compact.
Let $\eval\colon C(X,F)\times X\to F$,
$(f,x)\mto f(x)$ be the evaluation map. 
Then $\eval(M\times K)$
is compact.
\end{la}
\begin{proof}
The restriction map $\rho\colon C(X,F)_{c.o.}
\to C(K,F)_{c.o.}$, $f\mto f|_K$
being continuous by~\cite[\S3.2\,(2)]{Eng},
the set $\rho(M)$ is compact
in $C(K,F)_{c.o.}$.
The evaluation map
$\ve\colon C(K,F)\times K\to F$, $(f,x)\mto f(x)$
is continuous by \cite[Theorem~3.4.2]{Eng}.
Hence also
$\eval(M\times K)
=\ve(\rho(M)\times K)$ is compact.
\end{proof}
%
%
\begin{la}\label{hlp1}
Let $E$, $F_1$, $F_2$ and~$G$
be locally convex spaces,
and $\beta\colon F_1\times F_2\to G$
be a bilinear map which is hypocontinuous
with respect to compact subsets of~$F_2$.
Then
\[
C^n_\K(U,\beta )\colon C^n_\K(U,F_1)\times C^n_\K(U,F_2)\to
C^n_\K(U,G)\,,\quad (f,g)\mto \beta\circ (f, g)
\]
is hypocontinuous with respect
to compact subsets of $C^n_\K(U,F_2)$,
for each $n\in \N_0\cup\{\infty\}$.
If $\beta$ is continuous,
then also $C^n_\K(U,\beta)$ is continuous.
\end{la}
\begin{proof}
If $\beta$ is continuous
and hence smooth, then $C^n(U,\beta)$
is smooth and hence continuous,
as a very special case of
\cite[Proposition~4.16]{ZOO}.\footnote{Note that
the ordinary $C^n$-topology is used there, by
\cite[Proposition 4.19\,(d) and Lemma~A2]{ZOO}.}\\[2.5mm]
If $\beta$ is hypocontinuous,
it suffices to prove hypocontinuity of
$C^n(U,\beta)$ for each finite~$n$.
To see this, let $i_{n,G}\colon \colon C^\infty(U,G)\to C^n(U,G)$
be the inclusion map for $n\in \N_0$,
and define $i_{n,F_1}$ and $i_{n,F_2}$
analogously.
Since the topology on
$C^\infty(U,G)$ is initial
with respect to the maps $i_{n,G}$,
the restriction
$C^\infty(U,\beta)|_Y$
to a subset $Y\sub C^\infty(U,F_1)\times C^\infty(U,F_2)$
is continuous
if and only if
$i_{n,G}\circ C^\infty(U,\beta)|_Y
=C^n(U,\beta)|_{\wb{Y}}$
is continuous, where
$\wb{Y}:=(i_{n,F_1}\times i_{n,F_2})(Y)$.
Applying this with $Y=\{f\}\times C^\infty(U,F_2)$,
$Y=C^\infty(U,F_1)\times \{g\}$
and $Y=C^\infty(U,F_1)\times M$
with $M\sub C^\infty(U,F_2)$ compact,
we see that hypocontinuity
of each $C^n(U,\beta)$ implies
hypocontinuity of $C^\infty(U,\beta)$
(using the characterization given in
Proposition~\ref{prehypo}\,(c)).\\[2.5mm]
The case $n=0$.
Let $M\sub C(U,F_2)$ be compact
and consider a typical
$0$-neighbourhood in $C(U,G)$,
say $\lfloor K, W\rfloor$
with $K\sub U$ compact and
a $0$-neighbourhood $W\sub G$.
By Lemma~\ref{wannerf},
the set $N:=\eval(M\times K)\sub F_2$
is compact, where $\eval\colon C(U,F_2)\times U\to F_2$
is the evaluation map.
Since~$\beta$ is hypocontinuous,
there exists a $0$-neighbourhood
$V\sub F_1$ with
$\beta(V\times N )\sub W$.
Then
$\beta \circ (\lfloor K, V \rfloor\times M) \sub
\lfloor K, W\rfloor$
and hence $C(U,\beta)$ is hypocontinuous.\\[2.5mm]
Induction step. Given $n\in \N$,
assume that $C^{n-1}(U,\beta)$
is hypocontinuous in the second argument
for each~$U$ and~$\beta$.
The topology on $C^n(U,G)$ being initial with respect to the linear
maps $\lambda_1\colon C^n(U,G)\to C(U,G)_{c.o.}$,
$f \mto f$
and $\lambda_2\colon C^n(U,G)\to C^{n-1}(U\times E,G)$,
$f\mto df$
(by \cite[Lemma~A.1\,(d)]{ZOO}),
we only need to show that
$\lambda_j\circ C^n(U,\beta)$
is hypocontinuous
for $j\in \{1,2\}$.
We have $\lambda_1\circ C^n(U,\beta)=C(U,\beta)\circ (i_1\times i_2)$,
where $i_j\colon C^n(U,F_j) \to C(U,F_j)$
is the inclusion map which is continuous
and linear. Since $C(U,\beta)$ is hypocontinuous
by the case $n=0$,
we readily deduce that $\lambda_1\circ C^n(U,\beta)$
is hypocontinuous.
The map $\delta_j \colon C^n(U,F_j)\to C^{n-1}(U\times E,F_j)$,
$f\mto df$
is continuous linear
and $\pi\colon U\times E\to U$, $(x,y)\mto x$
is smooth, whence
$\rho_j \colon C^n(U,F_j)\to C^{n-1}(U\times E,F_j)$,
$\rho_j(f):=f\circ \pi$ is continuous
linear (cf.\ \cite[Lemma~4.4]{ZOO}).
By (\ref{abstrprod}),
we have
%
\begin{equation}\label{git}
\lambda_2\circ C^n(U,\beta)\;=\;
C^{n-1}(U\times E,\beta)\circ (\delta_1\times \rho_2)
\,+\, C^{n-1}(U\times E,\beta)\circ (\rho_1\times \delta_2)\,.
\end{equation}
Since $C^{n-1}(U\times E,\beta)\colon C^{n-1}(U\times E,F_1)
\times C^{n-1}(U\times E,F_2)\to C^{n-1}(U\times E,G)$
is hypocontinuous by the inductive hypothesis
and each summand in~(\ref{git})
is a composition thereof with a direct product
of continuous linear maps,
we deduce
that each summand and hence also $\lambda_2\circ C^n(U,\beta)$
is hypocontinuous in the
second argument. This completes the proof.
\end{proof}
%
%
\begin{la}\label{hlp2}
Let $E$, $F$ and~$G$
be locally convex spaces,
$U\sub E$ be open and
$\beta\colon E \times F \to G$
be a bilinear map which is hypocontinuous
with respect to compact subsets
of~$E$.
Then
\[
\beta_*\colon C^n_\K(U,F)
\to C^n_\K(U,G)\,,\quad
(\beta_* (f))(x)\, :=\, \beta(x,f(x))
\quad\mbox{for $\,f \in C^n_\K(U,F)$, $x\in U$}
\]
is a continuous linear map, for each $n\in \N_0\cup\{\infty\}$.
\end{la}
\begin{proof}
It is clear that~$\beta_*$ is linear.
We only need
to prove the assertion for finite~$n$,
by an argument similar to that in
the proof
of Lemma~\ref{hlp1}.
The proof is by induction.

The case $n=0$.
Suppose we are given a $0$-neighbourhood
in $C(U,G)$, say $\lfloor K,W\rfloor$
with $K\sub U$ compact and a $0$-neighbourhood
$W\sub G$. By Proposition~\ref{prehypo}\,(a),
there exists a $0$-neighbourhood
$V\sub F$ such that $\beta(K\times V)\sub W$.
Then $\beta_*(\lfloor K,V\rfloor)\sub \lfloor K,W\rfloor$.
Thus $\beta_*$ is continuous at~$0$ and hence
continuous, being linear.

Induction step. Let $n\in \N$ and assume
that the assertion holds for $n-1$ in place of~$n$.
Let us write $\beta_{*,n}\colon
C^n(U,F)\to C^n(U,G)$, for added clarity.
The topology on $C^n(U,G)$ being initial with respect to the linear
maps $\lambda_1\colon C^n(U,G)\to C(U,G)_{c.o.}$,
$f \mto f$
and $\lambda_2\colon C^n(U,G)\to C^{n-1}(U\times E,G)$,
$f\mto df$,
we only need to show that
$\lambda_j\circ \beta_{*,n}$
is continuous
for $j\in \{1,2\}$.
We have $\lambda_1\circ \beta_{*,n}
=\beta_{*,0}\circ i$,
where $i \colon C^n(U,F) \to C(U,F)$
is the continuous linear
inclusion map and $\beta_{*,0}$ is continuous
by the above.
To tackle $\lambda_2\circ \beta_{*,n}$,
note that
\[
A\colon (E\times E)\times F\to G\,,\quad
A((x,y),z)\,:=\,\beta(y,z)\quad\mbox{and}
\]
\[
B\colon (E\times E)\times F\to G\,,\quad
B((x,y),z)\,:=\,\beta(x,z)\quad\quad\;
\]
are hypocontinuous
with respect to compact subsets
of~$E\times E$.
Let $\pi_1\colon U\times E\to U$
be the projection on the first
component.
Then the map
$p\colon C^n(U,F)\to C^{n-1}(U\times E,F)$,
$f \mto f \circ \pi_1$
is continuous linear (cf.\ \cite[Lemma~4.4]{ZOO}).
By (\ref{abstrprod}), we have
%
\begin{equation}\label{lasttm}
\lambda_2\circ \beta_{*,n}\, =\,
A_*\circ p \, +\, B_*\circ d\,,
\end{equation}
where $d\colon C^n(U,F)\to C^{n-1}(U\times E,F)$,
$f \mto df$ is continuous linear and also
the maps $A_*\colon C^{n-1}(U\times E,F)\to
C^{n-1}(U\times E,G)$
and $B_*\colon C^{n-1}(U\times E,F)\to
C^{n-1}(U\times E,G)$
are continuous linear, by
the inductive hypothesis.
Hence $\lambda_2\circ \beta_{*,n}$
is continuous linear. This completes the proof.
\end{proof}
{\bf Proof of Theorem~\ref{poiss2}.}
By Lemma~\ref{hlp0},
the mapping $D\colon C^\infty(U,\K)\to C^\infty(U, E'_c)$,\linebreak
$f\mto f'$ is continuous and linear.
By Lemma~\ref{hlp1},
the bilinear map
\[
C^\infty(U,[.,.])\colon C^\infty(U,E')\times
C^\infty(U,E')\to C^\infty(U,E')\,,\quad
(f,g)\mto (x\mto [f(x),g(x)])
\]
is hypocontinuous with respect to compact subsets
of the second factor;
and if $[.,.]$ is continuous,
then also $C^\infty(U,[.,.])$.
The evaluation map
$\beta \colon E\times E'_c\to \K$, $\beta(x,\lambda):=\lambda(x)$
is hypocontinuous with respect to compact
subsets of~$E$
by Proposition~\ref{resulteval}.
Hence $\beta_*\colon C^\infty(U,E'_c)\to C^\infty(U,\K)$,
$f \mto \beta_*(f)=\beta\circ (\id_U,f)$
is continuous linear by Lemma~\ref{hlp2}.
Since
\[
\{.,.\}\, =\, \beta_*\circ C^\infty(U,[.,.])
\circ (D\times D)
\]
by definition,
we see that $\{.,.\}$ is a composition of
continuous maps if~$[.,.]$ is continuous,
and hence continuous.
In the general case,
$\{.,.\}$ is a composition of a hypocontinuous
bilinear map and continuous linear maps
and hence hypocontinuous.\,\Punkt
%
%
%
%
%
\section{Continuity of the map
taking {\boldmath $f$} to {\boldmath $X_f$}}\label{nwsec2}
In this section, we show continuity of the
mapping which takes
a smooth function to the corresponding Hamiltonian vector
field, in the case $\cS=c$.
\begin{thm}\label{poiss3}
Let $(E,[.,.])$ be a locally
convex Poisson vector space with respect to
$\cS=c$.
Let $U\sub E$ be an open subset.
Then the map
\begin{equation}\label{thmap}
\Psi \colon C^\infty_\K(U,\K)\to C^\infty_\K(U,E)\,,
\quad f\mto X_f
\end{equation}
is continuous and linear.
\end{thm}
\begin{proof}
Let $\eta_E\colon E\to (E'_c)'_c$
be the evaluation homomorphism
and
$V:= \{A\in L(E'_c,E'_c)\colon$
$(\forall x\in E)\; \eta_E(x)\circ A\in
\eta_E(E)\}$.
Then~$V$ is a vector subspace of $L(E'_c,E'_c)$
and $\ad(E')\sub V$.
The composition map
$\Gamma\colon (E'_c)'_c\times L(E'_c,E'_c)_c\to
(E'_c)'_c$, $(\alpha ,A)\mto \alpha\circ A$
is hypocontinuous with respect to equicontinuous
subsets of $(E'_c)'_c$,
by Proposition~9 in \cite[Chapter~III, \S5, no.\,5]{BTV}.
If $K\sub E$ is compact, then the polar
$K^\circ$ is a $0$-neighbourhood in~$E'_c$,
entailing that $(K^\circ)^\circ\sub (E'_c)'$
is equicontinuous.
Hence $\eta_E$ takes
compact subsets of~$E$ to equicontinuous subsets
of $(E'_c)'$, and hence
\[
\beta\colon E \times V \to E\,,\quad
(x,A) \mto \eta_E^{-1}(\Gamma(\eta_E(x), A))
\]
is hypocontinuous
with respect to compact subsets of~$E$.
Using Lemma~\ref{hlp2},
we see that $\beta_*\colon
C^\infty(U,V)\to C^\infty(U,E)$
is continuous linear.
Also the map
$D\colon C^\infty(U,\K)\to C^\infty(U,E'_c)$,
$f\mto f'$ is continuous linear by Lemma~\ref{hlp0}.
Furthermore,
$\ad=[.,.]^\vee \colon E'_c\to L(E'_c,E'_c)_c$
is continuous linear since
$[.,.]$ is hypocontinuous
(see Proposition~\ref{prehypo}\,(b)),
entailing that also
\[
C^\infty(U,\ad)\colon C^\infty(U,E'_c)\to
C^\infty(U,L(E'_c,E'_c)_c)\,,\quad
f \mto \ad\circ \, f
\]
is continuous linear
(see, e.g., \cite[Lemma~4.13]{ZOO}).
Hence
$\Psi= \beta_*\circ C^\infty(U,\ad)\circ D$
is continuous and linear.
\end{proof}
{\footnotesize
{\bf Helge Gl\"{o}ckner}, TU Darmstadt, Fachbereich Mathematik AG~5,
Schlossgartenstr.\,7,\\
64289 Darmstadt, Germany.
\,E-Mail:
\,{\tt gloeckner@mathematik.tu-darmstadt.de}}

\begin{thebibliography}{99}
%
%
\bibitem{Aus} Au\ss{}enhofer, L.,
{\em Contributions to the duality theory of Abelian
topological groups and to the theory
of nuclear groups}, Diss.\ Math.\ {\bf 384},
1999.
%
%
\bibitem{Ban} Banaszczyk, W., ``Additive Subgroups
of Topological Vector Spaces,'' Springer, 1991.
%
%
\bibitem{BGN}
Bertram, W., H. Gl\"{o}ckner and K.-H. Neeb,
\emph{Differential calculus over general base fields and rings},
Expo.\ Math.\ {\bf 22} (2004),
213--282.
%
%
\bibitem{BaS} Bochnak, J. and J. Siciak,
\emph{Analytic functions in topological vector spaces},
Studia Math.\ {\bf 39} (1971), 77--112.
%
%
\bibitem{BTV} Bourbaki, N., ``Topological Vector Spaces,
Chapters 1-5,'' Springer-Verlag, 1987.
%
%
\bibitem{BaO} Bunke, U. and M. Olbrich,
{\em Group cohomology and the singularities
of the Selberg zeta function associated to a Kleinian
group}, Ann.\ Math.\ {\bf 149} (1999), 627--689.
%
%
\bibitem{BO2} Bunke, U. and M. Olbrich,
\emph{The spectrum of Kleinian manifolds},
J.\ Funct.\ Anal.\ {\bf 172} (2000), 76--164.
%
%
\bibitem{Eng}
Engelking, R., ``General Topology,''
Heldermann Verlag, Berlin, 1989.
%
%
\bibitem{Flo}
Floret, K.,
\emph{Lokalkonvexe Sequenzen mit kompakten
Abbildungen}, J. Reine Angew.\
Math.\ {\bf 247} (1971), 155--195.
%
%
\bibitem{RES} Gl\"{o}ckner, H., \emph{Infinite-dimensional Lie groups
without completeness restrictions},
pp.\,43--59 in: A. Strasburger et al.\ (Eds.),
``Geometry and Analysis on Finite- and Infinite-Dimensional
Lie Groups,'' Banach Center Publications {\bf 55},
Warszawa, 2002.
%
%
\bibitem{FAM}
Gl\"{o}ckner, H., \emph{Remarks on holomorphic
families of operators}, TU~Darmstadt Preprint
{\bf 2256}, December~2002.
%
\bibitem{ZOO}
Gl\"{o}ckner, H., \emph{Lie groups over non-discrete
topological fields}, preprint,
arXiv:math.FA/0408008.
%
%
\bibitem{FUN}
Gl\"{o}ckner, H.,
{\em Fundamentals of direct limit Lie theory},
Compos.\ Math.\ {\bf 141} (2005), 1551--1577.
%
%
\bibitem{DL3}
Gl\"{o}ckner, H.,
\emph{Direct limits of infinite-dimensional Lie groups
compared to direct limits in related categories},
to appear in J. Funct.\ Anal.\
(cf.\ arXiv:math.FA/0606078).
%
%
\bibitem{COM}
Gl\"{o}ckner, H., \emph{Instructive
examples of smooth, complex differentiable
and complex analytic mappings
into locally convex spaces},
preprint, arXiv:math.FA/0701197.
%
%
\bibitem{HYP}
Gl\"{o}ckner, H.,
\emph{Bundles of locally convex spaces, group
actions, and
hypocontinuous bilinear mappings},
manuscript, Darmstadt 2002
(currently undergoing revision).
%
%
\bibitem{GGH}
Gl\"{o}ckner, H.,
R. Gramlich
and T. Hartnick,
\emph{Final group topologies, Phan systems
and Pontryagin duality},
preprint, cf.\ arXiv:math.GR/0603537.
%
%
\bibitem{GLN}
Gl\"{o}ckner, H., R.\,L. Lovas and K.-H. -Neeb,
\emph{Locally convex Poisson vector spaces
and co\-adjoint orbits} (provisional title),
work in progress.
%
%
\bibitem{GaN} Gl\"{o}ckner, H. and K.-H. -Neeb,
``Infinite-Dimensional Lie Groups. Vol.\,I:
Basic Theory and Main Examples,''
to appear 2007 or 2008 in Springer Verlag.
%
%
\bibitem{Gro}
Gro\ss{}e-Erdmann, K.-G.,
``The Borel-Okada Theorem Revisited,''
Habilitations\-schrift, FernUniversit\"{a}t
Hagen, 1992.
%
%
\bibitem{Jar}
Jarchow, H., ``Locally Convex Spaces,''
B.\,G. Teubner, Stuttgart, 1981.
%
%
\bibitem{Koe}
K\"{o}the, G., ``Topological Vector Spaces II,''
Springer-Verlag, New York, 1979.
%
%
\bibitem{KaM} Kriegl, A. and P.\,W. Michor,
``The Convenient Setting of Global Ana\-lysis,''
AMS, Providence, 1997.
%
%
\bibitem{Lan} Lang, S., ``Fundamentals of Differential
Geometry,'' Springer-Verlag, 1999.
%
%
\bibitem{Mai} Maissen, B.,
\emph{\"{U}ber Topologien im Endomorphismenraum
eines topologischen Vektorraumes},
Math.\ Ann.\ {\bf 151} (1963), 283--285.
%
%
\bibitem{Mil} Milnor, J.,
{\em Remarks on infinite-dimensional Lie groups},
pp.\,1007--1057 in:  DeWitt, B. and R. Stora (Eds.),
``Relativit\'{e}, Groupes et Topologie II,'' 1984.
%
%
\bibitem{Ne1}
Neeb, K.-H., \emph{Towards a Lie theory
of locally convex groups},
Japan J. Math.\ {\bf 1} (2006), 291--468.
%
%
\bibitem{Ne2} Neeb, K.-H.,
\emph{Poisson structures on pro-Lie groups},
manuscript in preparation.
%
%
\bibitem{OR1}
Odzijewicz, A. and T.\,S. Ratiu,
{\em Banach Lie-Poisson spaces and reduction},
Comm.\ Math.\ Phys.\ {\bf 243} (2003), 1--54.
%
%
\bibitem{OR2}
Odzijewicz, A. and T.\,S. Ratiu,
{\em Extensions of Banach Lie-Poisson spaces},
J. Funct.\ Anal.\ {\bf 217}
(2004), 103--125.
%
%
\bibitem{Sei}
Seip, U., ``Kompakt erzeugte Vektorr\"{a}ume
und Analysis,'' Springer Lecture Notes in Math.\
{\bf 273}, Springer, Berlin,1972.
%
%
\bibitem{Tho}
Thomas, E.\,G.\,F., ``Calculus on Locally Convex
Spaces,''
Preprint {\bf W-9604}
(unpublished),
Department of Mathematics,
University of Groningen, 1996.
%
%
\bibitem{Tre}
Treves, F., ``Topological Vector Spaces,
Distributions and Kernels,'' Academic Press, New York,
1967.
%
%
\bibitem{Wur}
Wurzbacher, T., \emph{Fermionic second quantization and the
geometry of the restricted Grassmannian},
pp.\,287--375 in:
A. Huckleberry and T. Wurzbacher (Eds.),
``Infinite-Dimensional K\"{a}hler
Manifolds,'' Birkh\"{a}user, Basel, 2001.
%
%
\end{thebibliography}
\end{document}